\documentclass[11pt]{pjm}
\input amssym.def
\input amssym.tex
\newcommand\RR{{\Bbb R}}

\def\d={\,:=\,}

\newcommand{\lpraum}[2]
 {{{\rm L}}^{#1}(#2)}

\font\frakten=eufm10
\newfam\frakfam
\textfont\frakfam=\frakten
\def\frk#1{\fam\frakfam{#1}}
\newtheorem{thm}{Theorem}
\newtheorem{lemma}[thm]{Lemma}

\newtheorem{Defn}[thm]{Definition}
\newtheorem{Ex}[thm]{Example}
\newtheorem{Rem}[thm]{Remark}
\newtheorem{Exs}[thm]{Examples}
\newtheorem{Rems}[thm]{Remarks}
\newtheorem{Defrem}[thm]{Definition and Remark}
\newtheorem{Remnt}[thm]{}

\newenvironment{rem}
 {\begin{Rem} \begin{rm}} {\end{rm} \hfill $\Box$ \end{Rem}}

\newenvironment{prf} {{\bf Proof.}}{\hfill $\Box$}

\begin{document}

\title{Hausdorff-Young inequalities for nonunimodular groups}
\author[H. F\"uhr]{Hartmut F\"uhr}
\address{ Institute of Biomathematics and Biometry \\
 GSF Research Center for Environment and Health \\D--85764 Neuherberg}
 \email{fuehr@mathematik.tu-muenchen.de}
 \urladdr{http://www.gsf.de/ibb/homepages/fuehr/}
\date{\today}
\begin{abstract} The paper studies Hausdorff-Young inequalities for certain group extensions,
by use of Mackey's theory. We consider the case in which the dual 
action of the quotient group is free almost everywhere. This
result applies to yield a Hausdorff-Young inequality for
nonunimodular groups.
\end{abstract}
\maketitle

\setcounter{section}{-1}
\section{Introduction}

In this paper we deal with a definition of ${\rm L}^p$-Fourier
transform on locally compact groups.
Recall that, for locally compact abelian groups, the Haus\-dorff-Young
inequality reads:\\
{\em Let $1 < p < 2$ and $q=p/(p-1)$. If $g \in {\rm L}^1(G)
 \cap {\rm L}^p(G)$, then $\widehat{g} \in {\rm L}^{q}(\widehat{G})$,
 with $\| \widehat{g} \|_{q} \le \| g \|_p$. }\\
The inequality allows to extend the Fourier transform to a continuous operator $
{\mathcal F}^p : {\rm
L}^p \to {\rm L}^q $ by continuity. It was generalized to type I unimodular groups by Kunze
\cite{Ku}. Over the years, various authors derived Hausdorff-Young
inequalities both for concrete groups \cite{KuSt,Ba,Fo,EyTe} and
for certain classes of groups \cite{Ru,Ru2,Ru3,In,BaLu}, with the
aim of getting a more precise bound in the inequality.

The formulation of the results for nonabelian groups requires a
certain amount of notation. Given a locally compact
group $G$, we denote by $\widehat{G}$ its unitary dual, i.e., the
set of (equivalence classes of) irreducible unitary representations,
endowed with the Mackey Borel structure. 
The dual space is used to define
the operator valued Fourier transform by letting
\[ {\rm L}^1 (G) \ni g \mapsto {\mathcal F}^1(g) := (\sigma(g))_{\sigma
\in \widehat{G}} ,\]
where $\sigma(g)$ is defined by the weak operator integral
\[ \sigma(g) = \int_G g(x) \sigma(x) dx ~~.\]
The main task for Plancherel theory in the nonabelian setting is
to find suitable Banach spaces for these operators, which is
what we describe next.

 For a Hilbert space ${\mathcal H}$, we denote
by ${\mathcal B}_p({\mathcal H})$ the set of bounded operators $T$ such
that $(T^*T)^{p/2}$ is trace class. ${\mathcal B}_p({\mathcal H})$ is a
Banach space with the obvious norm $\| T \|_p = {\rm
tr}((T^*T)^{p/2})^{1/p}$. The Plancherel theorem for type I unimodular locally
compact groups \cite[18.8.2]{Di} provides the
existence of a {\bf Plancherel measure} $\nu_G$ on $\widehat{G}$
with the following properties: Given a measurable realization
$(\sigma, {\mathcal H}_{\sigma})$ of representatives from
$\widehat{G}$, denote by ${\mathcal B}_p^{\oplus}$ ($p>1$) the space
of operator fields $A = (A_{\sigma})_{\sigma \in \widehat{G}}$
such that $A_{\sigma} \in {\mathcal B}_p({\mathcal H}_{\sigma})$ a.e., and
moreover that
\[ \left\| A \right\|_p := \left\| (A_\sigma)_{\sigma}
\right\|_{{\mathcal B}_p^{\oplus}} :=
\left( \int_{\widehat{G}} \left\| A_{\sigma} \right\|_p^p d\nu_G(\sigma)
\right)^{1/p} \]
is finite. It is routine to check that $({\mathcal B}_p^{\oplus},\| \cdot \|
_{{\mathcal B}_p^{\oplus}})$ is a Banach space; for $p=2$ even a Hilbert space.
It is convenient to use the notation ${\mathcal B}_{\infty}^{\oplus}$ for
the space of uniformly bounded operator fields.
Note that since $\| \pi(g) \| \le \| g \|_1$, for arbitrary representations
of $g$, we have the estimate
\begin{equation} \label{hdy1}
 \| {\mathcal F}(g) \|_{\infty} \le \| g \|_1
\end{equation}
Now, the Plancherel theorem for unimodular groups states that for the
unique Plancherel measure $\nu_G$ and for all
$g \in {\rm L}^2(G) \cap {\rm L}^1(G)$, we have
\begin{equation} \label{hdy2}
\left\| {\mathcal F}(g) \right\|_{{\mathcal B}_2^{\oplus}} = \| g \|_2 ~~.
\end{equation}
This gives rise to the {\bf Plancherel transform}
\[ {\mathcal P} = {\mathcal F}^2 : {\rm L}^2(G) \to {\mathcal B}_2^{\oplus} ~~,\]
which is a unitary equivalence.
The Hausdorff-Young inequality then states that, for all $1 \le
p \le 2$ and $g \in {\rm L}^p(G) \cap {\rm L}^1(G)$, we have
\begin{equation}
\label{hdy3}
 \left\| {\mathcal F}(g) \right\|_{{\mathcal B}_{q}^{\oplus}} \le \| g \|_p ~~,
\end{equation}
where $q=p/(p-1)$, which uniquely defines a Fourier transform 
${\mathcal F}^p : {\rm L}^p(G) \to {\mathcal B}_q^\oplus$. 
The proof of the inequality usually involves
interpolation techniques to derive the estimate from the
inequalities at the ``endpoints'', i.e., from (\ref{hdy1}) and
(\ref{hdy2}). However, for a sharper estimate of the operator norm 
of ${\mathcal F}_p$, in the following
denoted by $A_p(G)$, other techniques are required. 

For nonunimodular groups, an additional complication arises.
As Khalil showed for the {\bf ax+b}-group \cite{Kh}, there
exist $g \in {\mathcal C}_c(G)$ such that, for $\nu_G$-almost all
$\sigma$, ${\mathcal F}^1(g)(\sigma)$ is not a compact operator,
and hence ${\mathcal F}^1(g) \not\in {\mathcal B}_q^{\oplus}$, for
arbitrary $q$!

However, there is a way to save the Plancherel theorem, first obtained in
\cite{Kh} for the {\bf ax+b}-group and proved later on in full generality
in \cite{DuMo}. The Plancherel theorem in \cite{DuMo} can be
phrased as follows: Let $G$ be a nonunimodular locally compact group
with type I regular representation. There exists a field of unbounded
selfadjoint operators $(K_{\sigma})_{\sigma \in \widehat{G}}$, called
{\bf formal dimension operators}, and a measure $\nu_G$ on $\widehat{G}$,
such that, for every $g \in {\rm L}^2(G) \cap {\rm L}^1(G)$ and
$\nu_G$-almost every $\sigma$, the densely defined
operator $\sigma(g) K_{\sigma}^{1/2}$ extends to a Hilbert-Schmidt
operator, denoted $[\sigma(g) K_{\sigma}^{1/2}]$, and moreover
\[ \int_{\widehat{G}} \left\| [\sigma(g) K_{\sigma}^{1/2}] \right\|_2^2
 d\nu_G(\sigma) = \left\| g \right\|_2^2 ~~.\]
The operators
are unique up to scalar multiples, and once they are fixed, so is
the measure $\nu_G$. Similar results were obtained by \cite{KlLi,Ta}.

There does not seem to be a general treatment of the Hausdorff-Young
inequality for nonunimodular
groups. Apparently Terp proved such a result in a preprint
mentioned in \cite{In,BaLu}; however the paper seems never to have
appeared, and I have not been able to locate a copy of the
preprint. Several authors established Hausdorff-Young inequalities
for the {\bf ax+b}-group \cite{EyTe,Ru}. To my knowledge, the
largest class of groups under consideration were the solvable Lie
groups acting on Siegel domains, as investigated by Inoue
\cite{In}. The result is fairly intuitive, once we have seen how
to treat the ${\rm L}^2$-case: Letting
\[ {\mathcal F}^p(g)(\sigma) := [ \sigma(g) K_\sigma^{1/q}] \]
gives a well-defined operator from ${\rm L}^1(G) \cap {\rm
L}^p(G)$ into ${\mathcal B}_{q}^\oplus$ fulfilling
\[ \left\| {\mathcal F}^p(g) \right\|_{q} \le \| g \|_p ~~.\]

The first theorem shows that this result is true for
more general nonunimodular groups:
\begin{thm}
\label{thm:hdy_nonunimod}
 Let $G$ be a nonunimodular locally compact group such that
  $\lambda_G$ is type I and $N = {\rm Ker}(\Delta_G)$ is type I.
 Let $(K_{\sigma})_{\sigma \in \widehat{G}}$ denote the
 field of formal degree operators, and $\nu_G$ the Plancherel
 measure of $G$ belonging to that field. Let $1 < p < 2$ and
 $q = p/(p-1)$.

 Then, for all $g \in {\rm L}^1(G) \cap {\rm L}^p(G)$ and
 $\nu_G$-almost all $\sigma \in \widehat{G}$, the operator
 $\sigma(g) K_{\sigma}^{1/q}$ has a bounded extension
 $[ \sigma(g) K_{\sigma}^{1/q} ] \in {\mathcal B}_{q}({\mathcal H}_{\sigma})$, and
 we have the inequality
 \[ \left( \int_{\hat{G}} \|[ \sigma(g) K_{\sigma}^{1/q}]
 \|_{{\mathcal B}_{q}({\mathcal H}
 _{\sigma})}^{q} d \nu_G(\sigma) \right)^{1/q} \le A_p(N) \| g \|_p ~~,\]
 i.e., $A_p(G) \le A_p(N)$.
\end{thm}

Theorem 1 is a special case of a more general, somewhat
technical result (Theorem \ref{thm:hdy_grp_ext_A} below),
which proves a Hausdorff-Young inequality for certain group extensions
\[ 1 \to N \to G \to H \to 1 \]
with the additional properties that $N$ and $H$ are unimodular and
the dual action of $H$ on $\widehat{N}$ is free $\nu_{N}$-almost
everywhere. 

\section{Group extensions with free dual actions}

Now let $G$ be a second countable locally compact
group. $\mu_G$ denotes left-invariant Haar measure, and
$\Delta_G$ denotes the modular function of $G$.
We assume that $G$ is an extension of the unimodular normal
subgroup $N$, and that $H = G/N$ is unimodular as well.
Moreover we assume that $N$ is type I and $G$ has a type I regular representation,
and denote the Plancherel measures by $\nu_N$ and $\nu_G$, respectively. 
It is well-known that $H$ acts on $\widehat{N}$ via conjugation, denoted by
$H \times \widehat{N} \ni (\gamma, \sigma) \mapsto \gamma . \sigma$.
Mackey's theory allows the computation of $\widehat{G}$ from
this action, at least under certain regularity conditions.
As observed in \cite{KlLi}, the Mackey machine also provides the means to compute
Plancherel measure, under suitable assumptions on $N$ and the dual action. 
Our main technical assumption is the following:
\begin{enumerate}
 \item[(A)]
 There exists a Borel $\nu_N$-conull subset $U \subset \widehat{N}$ with the
 following property: $U$ is $H$-invariant with $U/H$ standard.
 Moreover, for all $\sigma_0 \in U$, ${\rm Ind}_N^G \sigma_0
 \in \widehat{G}$.
\end{enumerate}

In view of Mackey's theory, we obtain the following consequences:
\begin{enumerate}
\item[(i)] ${\rm Ind}_N^G \sigma_0
 \in \widehat{G}$ for every $\sigma_0 \in U$ implies that
 $H$ operates freely on the orbit $H . \sigma_0$. Moreover, ${\rm Ind}_N^G$ induces an
 injective map $\overline{{\rm Ind}}:U/H \hookrightarrow \widehat{G}$.
 \item[(ii)] Since $U/H$ is standard, there exists a measure decomposition along
 $H$-orbits (see for instance \cite[Theorem 2.1]{KlLi}).
 More precisely, there exist measures $\beta_{H. \sigma_0}$ on the orbits, such that 
 the Plancherel measures of $N$ and $G$ obey the relation
 \begin{equation} \label{eqn:disint_dual_0}
 d\nu_{N}(\rho) = d\beta_{H.\rho} (\rho)  d\nu_G(H.\rho)
 ~~.
 \end{equation}
 That the Plancherel measure of $G$ may be obtained as a quotient measure of 
 $\nu_N$, follows from
 Theorem \cite[I, 10.2]{KlLi}. For this we need to check four conditions:
 \begin{enumerate}
 \item[(I)] {\em $\lambda_N$ is type I with $\nu_N$ concentrated in $\widehat{N}_t$.} 
 This follows from the type I property of $N$.
 \item[(II)] {\em $\nu_N/H$ is countably separated.} This holds by Assumption (A).
 \item[(III)] {\em $\nu_N$-almost every little fixed group is trivial.} This holds by 
 Assumption (A).
 \item[(IV)] {\em $\lambda_G$ is type I.} This holds by assumption.
 \end{enumerate}
\item[(iii)] The fact that $H$ acts freely on $U$ entails that $U$
 can be seen as a product space $ H \times U_0$, by the following arguments:
 Since both $U$ and $U/H$ are standard, ($U$ as a Borel subset of $\widehat{N}$), 
 \cite[Theorem 5.2]{Ma} provides the existence of a Borel transversal, i.e. 
 a Borel subset $U_0 \subset U$ meeting each $H$-orbit in precisely one point.
 Hence the freeness of the 
 operation yields a measurable bijection $H \times U_0 \equiv U$, $(\gamma, \sigma_0) \mapsto \gamma
 . \sigma_0$, which by \cite[Theorem 3.2]{Ma} is a Borel isomorphism. 
 In the following, we use for $\sigma \in V$ the symbol $\sigma_0 \in U_0$
 to denote the representative of the associated dual orbit.
\item[(iv)] The Borel isomorphism $U \equiv H \times U_0$ allows to write down the
 measure disintegration (\ref{eqn:disint_dual_0}) much more explicitly.
 For this purpose we define the mapping $\psi: (h, \sigma_0)
 \mapsto \Delta_G(h)$ on $U$. Note that $\psi \equiv 1$ iff $G$ is unimodular.
 Then (\ref{eqn:disint_dual_0}) becomes
 \begin{equation} \label{eqn:disint_dual}
 d\nu_{N}(h, \sigma_0) =\psi(h) d\mu_{H}(h)  d\nu_G(\sigma)
 ~~.
 \end{equation}
 Indeed, a straightforward computation establishes
 that $\nu_{N}(\gamma . A) = \Delta_G(\gamma) \nu_N(A)$. This shows that 
 \[ d\beta_{H. \sigma_0}(\gamma . \sigma_0) = \alpha(H. \sigma_0) \psi(\gamma) d\mu_H(\gamma)~~,\]
 for some scalar factor $\alpha(H.\sigma_0)$. 
 Now $\nu_G$ can be renormalized to achieve that these factors are one.
\item[(v)]
 For the following calculations, $\sigma = {\rm Ind}_{N}^G \sigma_0$ is realized on
 ${\mathcal H}_\sigma = \lpraum{2}{H,$ $ d\mu_{H}; {\mathcal H}_{\sigma_0}}$ via cross sections
 (see Lemma \ref{lem:ind_cs}). 
 We define a family of operators
 $K_{\sigma}$ on ${\mathcal H}_\sigma$ given by multiplication with $\Delta_G$:
 \[ (K_{\sigma} \eta) (h) = \Delta_G(h) \eta(h) ~~,\]
 and ${\rm dom}(K_{\sigma})$ is the set of all $\eta \in \lpraum{2}{H,
 d\mu_{H};
 {\mathcal H}_{\sigma_0}}$, for which this product is also square integrable.

 Obviously $K_\sigma$ is the identity operator if $G$ is unimodular. In the other case,
 $K_\sigma$ is precisely the formal dimension operator, as can be seen by verifying
 the semi-invariance relation
 \[
 \sigma(x) K_\sigma \sigma(x)^* = \Delta_G(x)^{-1} K_\sigma~~,
 \]
 observing that the formal dimension operators obey the same relation \cite[Theorem 5]{DuMo}, and then
 applying the uniqueness statement \cite[Lemma 1]{DuMo}.
\end{enumerate}

\section{Hausdorff-Young inequalities for group extensions}

The following theorem is the main result of this paper.
\begin{thm}
\label{thm:hdy_grp_ext_A}
 Let a group extension $1\to N \to G \to H \to 1$ be given with $N,H$ unimodular,
 $N$ a type I group, and $\lambda_G$ a type I representation. 
 Assume that assumption (A) holds. Let $(K_{\sigma})_{\sigma \in \widehat{G}}$ denote the
 field of multiplication operators given in $(iv)$ above. Let $1 < p < 2$ and
 $q = p/(p-1)$.

 Then, for all $g \in {\rm L}^1(G) \cap {\rm L}^p(G)$ and
 $\nu_G$-almost all $\sigma \in \widehat{G}$, the operator
 $\sigma(g) K_{\sigma}^{1/q}$ has a bounded extension
 $[ \sigma(g) K_{\sigma}^{1/q} ] \in {\mathcal B}_{q}({\mathcal H}_{\sigma})$, and
 we have the inequality
 \[ \left( \int_{\hat{G}} \|[ \sigma(g) K_{\sigma}^{1/q}]
 \|_{{\mathcal B}_{q}({\mathcal H}
 _{\sigma})}^{q} d \nu_G(\sigma) \right)^{1/q} \le A_p(N) \| g \|_p ~~,\]
 i.e., $A_p(G) \le A_p(N)$.
\end{thm}

Before we prove this result, let us show how the nonunimodular
case follows from it:\\
{\bf Proof of Theorem \ref{thm:hdy_nonunimod}:} $N = {\rm
Ker}(\Delta_G)$ is a normal unimodular subgroup, with $H = G/N$
abelian. Assumption (A)
holds by \cite[Theorem 6]{DuMo}. Thus Theorem \ref{thm:hdy_grp_ext_A}
implies Theorem \ref{thm:hdy_nonunimod}. \hfill $\Box$.

We now proceed with the proof of Theorem \ref{thm:hdy_grp_ext_A}. It turns
out that it is a quite natural extension of the arguments used in \cite{Ru}
for the {\bf ax+b}-group, essentially by combining it with techniques from
\cite{KlLi}. First we show how $\sigma(f)$ acts via an operator
valued integral kernel. We then use an estimate of the $p$-norm of such operators
by certain cross norms, as provided by \cite{FoRu}. 

The first lemma computes the Haar measure of $G$ in terms of $\mu_N$ and $\mu_H$,
and fixes the normalizations we use in the following.
Note that, since $N$ is normal, $\Delta_G|_N = \Delta_N = 1$,
hence $\Delta_G$ can (and will) be regarded as a function on $H$.
\begin{lemma}
 Fix a measurable cross-section $\alpha: H \to G$. Then the mapping
 $N \times H \ni (n,h) \mapsto n \alpha(h) \in G$ is an isomorphism of Borel
 spaces.
 We use the notation $g = n \alpha(h) \equiv (n,h)$.
 Then
 \begin{equation} \label{eqn:HM_dec}
 d\mu_G(n,h) = d\mu_{N}(n) \Delta_G(h)
 d\mu_{H} (h)
\end{equation}
 is a left Haar measure.
\end{lemma}
\begin{prf} The map $(n,h) \mapsto n \alpha(h)$ is a measurable bijection between
 standard Borel spaces, and thus a Borel isomorphism \cite[Theorem 3.2]{Ma}.
 Fix $g = n \alpha(h), g' = n' \alpha(h') \in G$, then
 \[ g g' = n ~ \alpha(h) n' \alpha(h)^{-1} ~ \alpha(h) \alpha(h)' \alpha(h h')^{-1}~
 \alpha(h h') \]
 with $\alpha(h) n' \alpha(h)^{-1}, \alpha(h) \alpha(h)' \alpha(h h')^{-1} \in N$
 (observing $N \lhd G$). Hence right translation on $G$ corresponds to right translation in the
 variables $n,h$, though not by $n',h'$. Now the rightinvariance of $\mu_N,\mu_H$
 entails that $d\mu_N(n) d\mu_H (h)$ is a right Haar measure on $G$.
 But then $\Delta_G d\mu_N d\mu_H$ is a left Haar measure.
\end{prf}

The following two lemmas provide the integral kernels:

\begin{lemma} \label{lem:ind_cs}
 Let $\sigma_0 \in U_0$ and $\sigma = {\rm Ind}_{N}^G \sigma_0$. Define the cocycle
 $\Lambda: H \times H \to N$ by
 \[ \Lambda(\gamma,\xi) = \alpha(\xi)^{-1} \alpha(\gamma)
 \alpha(\alpha(\gamma)^{-1} \xi) ~~.\]
 If we realize $\sigma$ on
${\rm L}^2(H,d\mu_\gamma;{\mathcal H}_\sigma)$ via the cross-section $\alpha$, we obtain
 for $x = n \alpha(\gamma)$
\begin{equation} \label{eqn:real_ind_cs}
 \left( \sigma (x) f \right) (\xi) = (\sigma_0 \xi) (n) \sigma_0(\Lambda(\gamma,\xi))
 f(\gamma^{-1} \xi) ~~.\end{equation}
\end{lemma}

\begin{prf}
 Since the measure is invariant, the formula for induction via cross-sections yields
 \begin{eqnarray*}
  \left( \sigma (n,\gamma) f\right) (\xi)
  & = & \sigma_0 \left( \alpha(\xi)^{-1} n \alpha(\gamma) \alpha( \alpha(\gamma)^{-1}
  n^{-1} \xi ) \right) ~ f(\alpha(\gamma)^{-1} n^{-1} \xi) \\
 & = & \sigma_0 \left( \alpha(\xi)^{-1} n \alpha(\xi) \Lambda(\gamma,\xi) \right) ~ f(\gamma^{-1} \xi) \\
 & = &  (\xi . \sigma_0) (n) \sigma_0(\Lambda(\gamma,\xi))
 f(\gamma^{-1} \xi)~~,
 \end{eqnarray*}
 where $\alpha(\gamma)^{-1} n^{-1} \xi  = \alpha(\gamma)^{-1} \xi$ is
 due to $N \lhd G$ .
\end{prf}

The next step consists in integrating this representation:

\begin{lemma}
\label{lem:op_kernel}
 Let $ \sigma_0 \in U_0$ and $\sigma = {\rm Ind}_{N}^G \sigma_0$.
 For $g \in {\rm L}^1(G)$ let $g_\gamma := g(\cdot,\gamma)$.
 Then $\sigma(g): {\rm L}^2(H,d\mu_{H};{\mathcal H}_{\sigma_0})
 \to  {\rm L}^2(H,d\mu_{H};{\mathcal H}_{\sigma_0})$
 can be written as
 \[ \sigma(g) f (\xi) = \int_{H} k_\sigma(\xi,\gamma) f(\gamma)
 d \gamma ~~,\]
 where $k_\sigma$ is an operator valued integral kernel given by
 \[ k_\sigma(\xi,\gamma) = (\xi . \sigma_0) (g_{\xi \gamma^{-1}}) \circ
 \sigma_0(\Lambda(\xi \gamma^{-1},\xi)) \cdot \Delta_G(\xi \gamma^{-1}) ~~.\]
\end{lemma}

\begin{prf}
 First note that by Fubini's theorem
 $g_{\gamma} \in {\rm L}^1(N)$, for almost every $\gamma \in H$,
 which justifies the use of $ (\xi . \sigma_0) (g_{\xi \gamma^{-1}})$.
 The following formal calculations can be made rigorous by plugging
 them into scalar products, according to the definition of the
 weak operator integral. Using the previous lemma and unimodularity of $H$, we see that
 \begin{eqnarray*}
 (\sigma(g) f)(\xi) & = &
 \int_H \int_{N} g(n,\gamma) (\xi .  \sigma_0) (n) \sigma_0(\Lambda(\gamma,\xi))
 f(\gamma^{-1} \xi)
 dn \Delta_G(\gamma) d\gamma \\
 & = & \int_H (\xi . \sigma_0) (g_{\gamma}) \sigma_0(\Lambda(\gamma,\xi)) f(\gamma^{-1} \xi)
 \Delta_G(\gamma) d\gamma \\
 & = & \int_H (\xi . \sigma_0) (g_{\gamma^{-1}}) \sigma_0(\Lambda(\gamma^{-1},\xi)) f(\gamma \xi)
 \Delta_G(\gamma^{-1}) d\gamma \\
 & = &  \int_H (\xi . \sigma_0) (g_{\xi \gamma^{-1}}) \sigma_0(\Lambda(\xi \gamma^{-1},\xi))
 \Delta_G(\xi \gamma^{-1}) f(\gamma) d\gamma~~,
 \end{eqnarray*}
 which is the desired formula.
\end{prf}

{\bf Proof of Theorem \ref{thm:hdy_grp_ext_A}:}
 Let $g \in {\rm L}^p(G) \cap {\rm L}^1(G)$ be given. By Lemma
 \ref{lem:op_kernel} and the definition of $K_\sigma$, we find
 that $\sigma(g) K_{\sigma}^{1/q}$ has the operator-valued
 kernel
 \[ k_{\sigma}(\xi,\gamma)  = (\xi . \sigma_0) (g_{\xi \gamma^{-1}}) \sigma_0(\Lambda(\xi \gamma^{-1},\xi))
 \Delta_G(\xi \gamma^{-1}) \Delta_G(\gamma)^{1/q} ~~.\]
 We want to use a result from \cite{FoRu}, which gives
 an estimate of $ \| \sigma (g) K_{\sigma}^{1/q} \|_p$ in
 terms of the {\bf cross norm}
 \[ \left\| \tilde{k}_{\sigma} \right\|_{q,p,q} :=
 \left( \int_{H} \left[ \int_{H} \| k
 (\xi,\gamma) \|_{q}^p d\xi \right]^{q/p}  d \gamma \right)^{1/q} ~~,\]
 for arbitrary operator-valued kernels.
 Using first \cite[Corollary 1]{FoRu} and then the Cauchy-Schwarz-inequality, we have
 \begin{eqnarray} \nonumber
 \int_{\widehat{G}} \left\| \sigma(g) K_{\sigma}^{1/q} \right\|_{q}^{q}
 d\nu_G(\sigma) & \le &
 \int_{\widehat{G}} \left\| k_{\sigma} \right\|_{q,p,q}^{q/2}
 \left\| k_{\sigma}^* \right\|_{q,p,q}^{q/2}  d\nu_G(\sigma) \\
 & \le &  \left( \int_{\widehat{G}} \left\| k_{\sigma}
 \right\|_{q,p,q}^{q}  d\nu_G(\sigma) \right)^{1/2}
  \left( \int_{\widehat{G}} \left\| k_{\sigma}^*
 \right\|_{q,p,q}^{q}  d\nu_G(\sigma) \right)^{1/2}~~, \label{cn_est}
 \end{eqnarray}
 where $k_{\sigma}^*(\xi, \gamma) = k_{\sigma}(\gamma,\xi)^*$.
 It remains thus to estimate the integral over the cross norms.
 We have that
 \begin{eqnarray*}
 \lefteqn{ \int_{\widehat{G}} \left\| k_{\sigma}
 \right\|_{q,p,q}^{q}  d\nu_G(\sigma) =} \\ & = &
 \int_{U_0}  \int_{H} \left[ \int_{H}
 \| (\xi . \sigma_0) (g_{\xi \gamma^{-1}})  \sigma_0(\Lambda(\xi \gamma^{-1},\xi))
 \Delta_G(\xi \gamma^{-1})  \Delta_G(\gamma)^{1/q}
 \|_{q}^p d\xi \right]^{q/p}
 d\gamma d\nu_G(\sigma) \\
 & = &
 \int_{U_0}  \int_{H} \left[ \int_{H}
 \| (\xi \gamma .\sigma_0) (g_{\xi}) \Delta_G(\xi) \Delta_G(\gamma)^{1/q}
 \|_{q}^p d \xi \right]^{q/p}
 d\gamma d\nu_G(\sigma) \\
 & \le & \left( \int_{H} \left[ \int_{U_0} \int_{H}
 \| (\xi \gamma .\sigma_0) (g_{\xi}) \Delta_G(\xi) \Delta_G(\gamma)^{1/q}
 \|_{q}^{q} d\gamma d\nu_{G}(\sigma) \right]^{p/q}
 d \xi \right)^{q/p} \\
 & = & \left(\int_{H} \left[ \int_{U_0} \int_{H}
 \| (\gamma .\sigma_0) (g_{\xi}) \Delta_G(\xi) \Delta_G(\xi^{-1} \gamma)^{1/q}
 \|_{q}^{q} d\gamma d\nu_{G}(\sigma) \right]^{p/q}
 d\xi \right)^{q/p} \\
 & = &  \left(\int_{H} \left[ \int_{U_0} \int_{H}
 \| (\gamma .\sigma_0) (g_{\xi}) \|_{q}^{q} \Delta_G(\gamma)
 d\gamma d\nu_{G}(\sigma)
 \right]^{p/q}
 \Delta_G(\xi) d\xi \right)^{q/p}~~.
 \end{eqnarray*}
 Note that we have tacitly dropped the unitary operators
 $ \sigma_0(\Lambda(\xi \gamma^{-1},\xi))$, since they obviously
 do not affect the $\| \cdot \|_{q}$-Norm.
 The inequality is due to Minkowski's generalized inequality.
 Now, by the measure disintegration (\ref{eqn:disint_dual}), we find that
 the inner double integral can be estimated by use of the Hausdorff-Young
 inequality for $N$:
 \begin{eqnarray*}
 \lefteqn{\left(\int_{H} \left[ \int_{U_0} \int_{H}
 \| (\gamma . \sigma_0) (g_{\xi}) \|_{q}^{q} \Delta_G(\gamma)
 d\gamma d\nu_{G}(\sigma)
 \right]^{p/q} \Delta_G(\xi) d\xi \right)^{q/p}} \\ & = &
 \left( \int_{H} \| {\mathcal F}^p(g_{\xi}) \|_{q}^{p}
 \Delta_G(\xi) d\xi \right)^{q/p} \\
 & \le & \left(  \int_{H} A_p^{p} \| g_{\xi} \|_p^p
 \Delta_G(\xi) d\xi \right)^{q/p} \\
 & = & A_p(N)^{q} \| g \|_p^{q}~~,
 \end{eqnarray*}
 which takes care of the first factor in (\ref{cn_est}).
 For the second factor, we have to compute the cross norms of
 \[ k_{\sigma}^*(\xi,\gamma) =  \sigma_0(\Lambda( \gamma \xi^{-1},\gamma))^*
 \circ (\gamma . \sigma_0) (g_{\xi^{-1} \gamma})^* \Delta_G(\gamma \xi^{-1})
 \Delta_G(\xi)^{1/q} ~~.\]
 Here we see that
 \begin{eqnarray*}
  \lefteqn{\int_{\widehat{G}} \left\| k_{\sigma}^*
 \right\|_{q,p,q}^{q}  d\nu_G(\sigma) = } \\ & = &
 \int_{U_0}  \int_{H} \left[ \int_{H}
 \| \sigma_0(\Lambda(\gamma \xi^{-1},\gamma))^* (\gamma . \sigma_0) (g_{\xi^{-1} \gamma})^*
 \Delta_G(\xi^{-1} \gamma)  \Delta_G(\xi)^{1/q}
 \|_{q}^p d\xi \right]^{q/p}
 d\gamma d\nu_G(\sigma) \\
 & = &
 \int_{U_0}  \int_{H} \left[ \int_{H}
 \| (\gamma . \sigma_0) (g_{\xi^{-1}})^*
 \Delta_G(\xi^{-1})  \Delta_G(\xi \gamma)^{1/q}
 \|_{q}^p d\xi \right]^{q/p}
 d\gamma d\nu_G(\sigma) \\
 & \le & \left( \int_{H} \left[ \int_{U_0} \int_{H}
 \| (\gamma . \sigma_0) (g_{\xi^{-1}})^* \Delta_G(\xi)^{-1+1/q}
 \Delta_G(\gamma)^{1/q} \|_{q}^{q} d\gamma d\nu_G(\sigma)
 \right]^{p/q} d\xi \right)^{q/p} \\
 & \le &  \left( \int_{H} \left[ \int_{U_0} \int_{H}
 \| (\gamma . \sigma_0) (g_{\xi})^* \Delta_G(\xi)^{1-1/q}
 \Delta_G(\gamma)^{1/q} \|_{q}^{q} d\gamma d\nu_G(\sigma)
 \right]^{p/q} d\xi \right)^{q/p} \\
 & = &  \left(\int_{H} \left[ \int_{U_0} \int_{H}
 \| (\gamma .\sigma_0) (g_{\xi}) \|_{q}^{q} \Delta_G(\gamma)
 d\gamma d\nu_{G}(\sigma)
 \right]^{p/q}
 \Delta_G(\xi) d\xi \right)^{q/p}~~,
 \end{eqnarray*}
 where the inequality is again the generalized Minkowski inequality,
 and we have used that taking adjoints and multiplication with unitaries are
 isometries on ${\mathcal B}_p$. Now we can
 conclude in the same way as before. \hfill $\Box$

 \begin{rem} \label{rem:L2_case}
 Note that  for the case $p=2$, all inequalities are in fact equalities:
 That the Hilbert-Schmidt norm of an operator given by an ${\rm L}^2$-kernel
 equals the ${\rm L}^2$-norm of the kernel is well-known, i.e., instead of
 (\ref{cn_est}) we have
 \[ \int_{\widehat{G}} \| \sigma(f) K_{\sigma}^{1/2} \|_2^2 d\nu_G(\sigma) =
 \int_{\widehat{G}} \| k_{\sigma} \|_{2}^2 d\nu_G(\sigma) ~~.\]
 Instead of the generalized Minkowski inequality, we can simply
 apply Fubini's theorem (since $p=q=2$), replacing the ``$\le$'' by
 ``$=$''. The last inequality in the argument is now an instance of
 the Plancherel theorem for $N$, hence once more an equality.
 Hence the computation provides a rather concrete proof that $\overline{\nu}=\nu_G$, 
 once the two measures are proven to be equivalent.
\end{rem}

\begin{rem}
 The type I assumption on $N$ ensures that we may define the spaces ${\mathcal B}_q^\oplus$
 as direct integrals of Schatten-von-Neumann spaces. In the general case, other traces
 than the natural operator trace may occur. For this setting, an extension of the
 norm estimates for Hilbert-Schmidt-valued kernels by cross norms, as obtained in \cite{Ru},
 to more general traces will be required. 
\end{rem}


\begin{rem}
 Another situation where Theorem \ref{thm:hdy_grp_ext_A} is applicable occurs in the
 context of simply connected, connected nilpotent Lie groups.  Baklouti, Smaoui and
 Ludwig \cite{BaLu} proved for these groups the estimate
 \begin{equation} \label{eqn:hdy_nilp} A_p(G) \le A_p(\RR)^{{\rm dim}(G) - d^*(G)/2} ~~,\end{equation}
 where $d^*(G)$ is the maximal coadjoint orbit dimension. In particular, let 
 $N \lhd G$ be a connected, codimension 1 normal subgroup of the simply connected, connected
 Lie group $G$. Assume that $d^*(N)< d^*(G)$, i.e., $d^*(N)=d^*(G)-2$. Given  
 $l \in {\frk g}^*$ with maximal orbit dimension, let $l_0$ denote the restriction to ${\frk n}$, 
 and let
 \[ {\frk r}_l = \{ X \in {\frk g} : l([X,Y]) = 0 \forall Y \in {\frk g} \} \]
 denote the {\bf radical} of $l$ in ${\frk g}$. Denote by ${\frk r}_{l_0} \subset {\frk n}$
 the radical of $l_0$ in ${\frk n}$. Then we have that 
 \[ {\rm dim}({\frk r}_l) =
\dim (G) - d^*(G) \] and 
 \[ {\rm dim}({\frk r}_{l_0}) = \dim(N) - \dim({\mathcal O}_{l_0}) \ge \dim(G) - 1 - (d^*(G) - 2)
 = \dim(G) - d^*(G)+1 ~~.
 \]
  Thus ${\frk r}_l
 \subset {\frk n}$ by \cite[Proposition 1.3.4]{CoGr}, and 
 \cite[Theorem 2.5.1(a)]{CoGr} implies that $\pi_l = {\rm Ind}_N^G \pi_{l_0}$.
 Here $\pi_l \in \widehat{G}$ and $\pi_{l_0} \in \widehat{N}$ denote the
 representations associated to $l,l_0$ by Kirillov's construction. 
 This holds for all $l_0$ for which ${\rm dim}{\mathcal O}_l= d^*(G)$.
 But then the set
 \[
 U = \{ \pi_{l_0} : \dim ({\mathcal O}_l) = d^*(G) \} \subset \widehat{N}
 \]
 is conull and $G$-invariant. Here we used the notation $\pi_l$ for the representation associated
 to $l$ by the Kirillov construction. 
 Hence, assuming (\ref{eqn:hdy_nilp}) for $N$, we obtain from Theorem \ref{thm:hdy_grp_ext_A} that
 \begin{eqnarray*} A_p(G) &\le & A_p(N) = A_p(\RR)^{{\rm dim}(N) - d^*(N)/2} = A_p(\RR)^{{\rm dim}(G) - 1 - (d^*(G)-2)/2} \\
 & = & A_p(\RR)^{{\rm dim}(G)-d^*(G)/2} ~~.\end{eqnarray*}
 In other words, Theorem \ref{thm:hdy_grp_ext_A} provides one half of the induction step
 for the proof of (\ref{eqn:hdy_nilp}). The other half would have to deal with the case
 $d^*(N) = d^*(G)$, which by similar arguments as above implies that the dual action 
 of $G/N$ on $\widehat{N}$ is trivial almost everywhere. 
\end{rem}

\end{document}